%% file: main.tex
\documentclass{amsart}
\usepackage{amsmath,amssymb}
\usepackage{tikz} 
\usepackage{graphicx}
\newcommand*\Laplace{\mathop{}\!\mathbin\bigtriangleup}

\newtheorem{theorem}{Theorem}

\newtheorem{proposition}{Proposition}[section]

\newtheorem{lemma}[proposition]{Lemma}
\newtheorem*{main-theorem}{Main Theorem}
\newtheorem*{theorem*}{Theorem}

\theoremstyle{definition}

\newtheorem{remark}[proposition]{Remark}

\newtheorem*{remark*}{Remark}

\numberwithin{equation}{section}

\def\O{{\mathcal O}}
\def\reals{{\mathbb R}}
\def\p{\partial}
\begin{document}

\title{Neumann Data Mass on perturbed triangles}
%
%
\author[H. Christianson]{Hans Christianson}
\address[H. Christianson]{ Department of Mathematics, University of North Carolina.\medskip}
 \email{hans@math.unc.edu}

 \author[J. Xi]{Jin Xi}
\address[J. Xi]{ Department of Economics, University of California San
  Diego.\medskip}
 \email{x5jin@ucsd.edu}

%
%
%

\begin{abstract}
Based on a previous paper 
 \cite{Chr-tri} on Neumann data for Dirichlet eigenfunctions on
 triangles, we extend the study in two
ways. First, we investigate the (semi-classical) Neumann data mass on
perturbed triangles. Specifically, we replace one side of a triangle
by adding a smooth perturbation, and assume that the disparity between
the perturbation and the original side is bounded by a small value
$\epsilon$.  Second, we add a small $\epsilon$ sized potential to the (semi-classical)
Laplacian  and see how the results change on triangles. In both cases, we find that the $L^2$ norm of Neumann data on each side is close to the length of the side divided by the area of the triangle, and the difference is  dominated by $\epsilon$. 
\end{abstract}
\maketitle              

\section{Introduction}
The purpose of this paper is to extend the results on triangles in \cite{Chr-tri}
by the first author to small perturbations, both in the domain and in
the Laplacian.

In our first extension, we study the Dirichlet eigenfunction problem
in a new domain that is modified from a triangle.  Given a bounded domain $D
\subset \reals^2$, consider the Dirichlet eigenfunction problem: 
\begin{equation}
\begin{cases}
 &(-h^2\Laplace - 1)u = 0 \text{ in } D \\
 &u|_{\partial D} = 0
\end{cases}
\label{E:efcn}
\end{equation}
where the eigenfunctions are assumed to be normalized $\|u\|_{L^2 (D)} = 1$.

In \cite{Chr-tri} by the first author, it is  shown that if the domain $D$ is a triangle,  the Neumann
data mass is equally distributed on each face: If $F_1, F_2$, and
$F_3$ are the three faces of $T$ with lengths $l_1, l_2,$ and $l_3$
respectively, then 
\begin{equation*}
\int_{F_j}|h\partial_{\nu}u|^2dS = \frac{l_j}{\text{Area}(D)}, \quad j
= 1,2,3,
\end{equation*}
where $h\partial_{\nu}u$ is the semi-classical normal derivative on
$\partial T$, $dS$ is the arclength measure,  $\text{Area}(D)$ is the
area of the triangle $D$. 
 In other words, the $L^2$ norm of Neumann data on
each side equals to the length of the side divided by the area of the
triangle.  An analogous result holds \cite{Chr-simp} when the dimension $n \geq 3$, so
the purpose of the present paper is to extend the 2 dimensional
results to other domains.

To look into the Neumann data mass in different domains, we first
study domains that are close to triangles: we construct a  planar
domain by changing one side of a triangle to a smooth function which
is close to linear in a suitable sense. Moreover, we restrict this
function to be close to the original side of the triangle, and their
disparity is bounded by a small number $\epsilon$. Our main findings
show that the Neumann data mass on each side is close to that of the
original triangle, and the difference is dominated by
$\epsilon$.

\begin{figure}
\hfill
\centerline{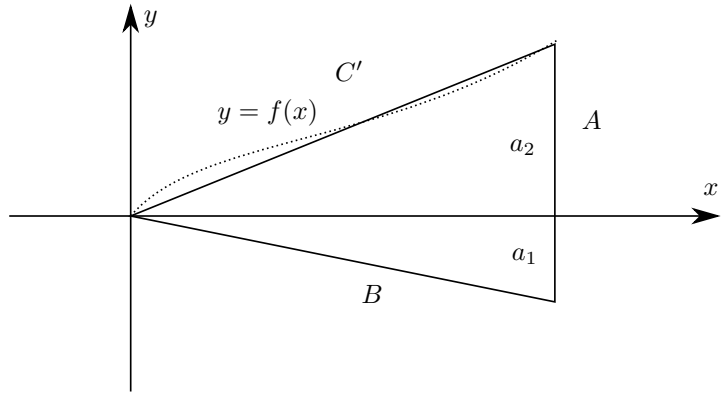}
\caption{Setup for acute (and right) triangles.}\label{exampleplot}
\hfill
\end{figure}

\begin{figure}
\hfill
\centerline{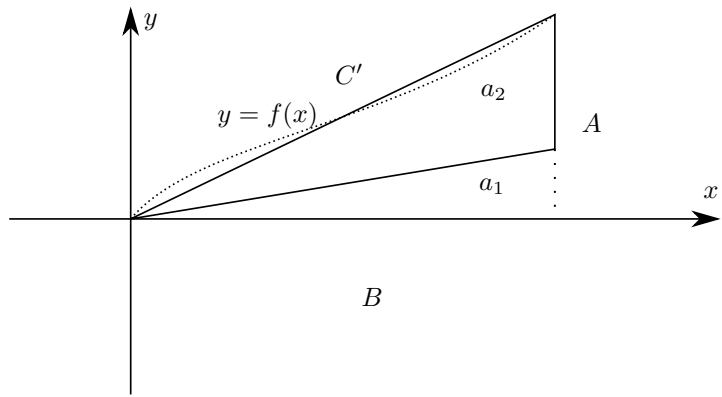}
\caption{Setup for obtuse triangles.}\label{exampleplot2}
\hfill
\end{figure}

Let $T \subset \reals^2$ be a  triangle with sides $A,$ $B,$ and
$C$ with lengths $a,$ $b$,  and $c$ respectively.  Assume the triangle
is oriented so that one corner is at the origin, side $A$ is vertical,
and side $C$ is parametrized by $a_2 x/l$, $0 \leq x \leq l$,
for positive $a_2$ and $l$ (see Figure \ref{exampleplot} for a picture in the
acute case and Figure \ref{exampleplot2} for the obtuse case).  Let $\tilde{g}(x)$ be a
smooth function satisfying $\tilde{g}(0) = \tilde{g}(l) = 0$, $| \tilde{g}(x) | \leq 1$, and $| \tilde{g}'(x) |
\leq 1$.   For $\epsilon >0$ small, let $g_\epsilon(x) = \epsilon \tilde{g}(x)$.  Let $D_\epsilon \subset \reals^2$
be the domain with side $C$ replaced by side $C'$ parametrized by $a_2
x/l + g_\epsilon(x)$.

\begin{theorem}
  \label{T:1}
  Fix $\epsilon>0$ small and 
let $D_\epsilon \subset \reals^2$ be the domain described above, and
suppose $\{u^\epsilon_h\}_h$ solves the semiclassical eigenfunction problem
\eqref{E:efcn} with $D$ replaced by $D_\epsilon$.  Then 
\begin{equation*}
\int_{A}|h\partial_{\nu}u^\epsilon |^2dS = \frac{a}{\text{Area}(D)}  +\O(\epsilon)
\end{equation*}
\begin{equation*}
\int_{B}|h\partial_{\nu}u^\epsilon|^2dS = \frac{b}{\text{Area}(D)}  + \O(\epsilon),
\end{equation*}
and
\begin{equation*}
\int_{C'}|h\partial_{\nu}u^\epsilon|^2dS = \frac{l(C')}{\text{Area}(D)} + O(\epsilon),
\end{equation*}
where $l(C')$ is the length of side $C'$.
\end{theorem}
\begin{remark}
  The theorem states that the equidistribution law from \cite{Chr-tri}
  is stable under small perturbations.  The implicit constants are
  independent of $\tilde{g}$ as long as it satisfies $| \tilde{g} |
  \leq 1$ and $| \tilde{g}' | \leq 1$.  In fact, we really only use
  that $| g | \leq \epsilon$ and $| g' | \leq \epsilon$ in the proof,
  so $g$ does not have to be of the form $\epsilon \tilde{g}$, however
  it does make the statement of the theorem more clear.

  \end{remark}

\begin{remark}
  We expect that the analogue of Theorem \ref{T:1} (and Theorem
  \ref{T:2} below) hold in higher dimensions, following the work
  \cite{Chr-simp} by the first author.
  \end{remark}

In our second extension, we consider a modified Dirichlet
eigenfunction problem on triangles.
\begin{theorem}
  \label{T:2}
  Let $T \subset \reals^2$ be a triangle with sides $A$, $B$, and
  $C$ with lengths $a$, $b$, and $c$ respectively.
Let $\tilde{w}(x,y)$ be a smooth function on $T$ with $| \tilde{w} |
\leq 1$ and $| \nabla \tilde{w} | \leq 1$.  For $\epsilon>0$ small,
let $w_\epsilon = \epsilon \tilde{w}$.  
  Consider the eigenfunction problem
\begin{equation*}
\begin{cases} 	-h^2\Laplace + w_\epsilon(x,y))u^\epsilon = u^\epsilon \text{ on } T,\\
 	u^\epsilon|_{\partial_T} = 0,
\end{cases}
\end{equation*}
and assume the $u^\epsilon$ are normalized $\| u^\epsilon\|_{L^2(T)}
=1$.  
Then for $\epsilon>0$ sufficiently small
\begin{equation*}
\int_{A}|h\partial_{\nu}u^\epsilon |^2dS = \frac{a}{\text{Area}(T)}  +\O(\epsilon)
\end{equation*}
\begin{equation*}
\int_{B}|h\partial_{\nu}u^\epsilon|^2dS = \frac{b}{\text{Area}(T)}  + \O(\epsilon),
\end{equation*}
and
\begin{equation*}
\int_{C}|h\partial_{\nu}u^\epsilon|^2dS = \frac{c}{\text{Area}(T)} + O(\epsilon).
\end{equation*}
\end{theorem}

%

\subsection{History}
The study of restrictions of eigenfunctions and the study of boundary
traces is an old subject.  In this very abbreviated history we just
focus on some of the  recent developments particularly relevant to the
present work.
  Previous results on restrictions primarily focused on upper bounds.
  In general, it is difficult to separate the behaviour of the
  Dirichlet and Neumann data for restrictions to interior
  hypersurfaces.  
  In the paper of Burq-G\'erard-Tzvetkov \cite{BGT-erest}, restrictions of
  the Dirichlet data to arbitrary smooth hypersurfaces on manifolds
  without boundary were considered.  An upper
  bound of the norm (squared) of the restrictions of $\O(h^{-1/2})$ was
  proved, and shown to be sharp.  
  Of course this sharpness shows that there are
  {\it some} eigenfunctions with a known lower bound.   In
  the first  author's paper with Hassell-Toth \cite{CHT-ND}, an upper bound of $\O(1)$ was
  proved for (semi-classical) Neumann data restricted to arbitrary
  smooth hypersurfaces on manifolds without boundary,  and this also shown to be sharp.  Again, this gives a lower and
  upper bound for the Neumann data alone for {\it some} eigenfunctions.

  In the case of quantum ergodic eigenfunctions,  more is
  known.  In the papers of G\'erard-Leichtnam \cite{GeLe-qe} and Hassell-Zelditch \cite{HaZe}, the Neumann (respectively
  Dirichlet) boundary data of Dirichlet (respectively Neumann) quantum ergodic eigenfunctions is
  studied, and shown to have an asymptotic formula for a {\it density
    one} subsequence.  
  Similar statements were proved for
  interior hypersurfaces by Toth-Zelditch \cite{ToZe-1, ToZe-2}.
  Again, potentially a sparse subsequence may behave differently.  In
  the  author's paper with Toth-Zelditch \cite{CTZ-1}, an asymptotic formula for the
  whole weighted Cauchy data is proved for the entire sequence of
  quantum ergodic eigenfunctions, however it is impossible to separate
  the behaviour of the Dirichlet versus Neumann data.  And the
  sequence of eigenfunctions is assumed to already be quantum ergodic,
  having thrown out any weird sparse subsequences.

  In \cite{Chr-tri} and \cite{Chr-simp}, the first author studied the
  Neumann boundary data for Dirichlet eigenfunctions on simplices and proved an
  equidistribution law.  It is not an asymptotic, but an exact
  identity, and holds for the entire sequence of eigenfunctions.  It
  agrees with what the paper of Hassell-Zelditch would give, but is an
  exact identity and holds for the whole sequence of eigenfunctions, so hints
  at quantum ergodicity (or at least some weak form of quantum
  ergodicity).  The main purpose of this paper is
to study similar phenomena for eigenfunction problems which are
``close to'' triangle eigenfunction problems.

\section{Proof of Theorem \ref{T:1}}

For the rest of the paper, let us drop the awkward $u^\epsilon$
notation and just write $u$, being careful to keep in mind that
everything implicitly depends on $\epsilon$.  

We first need to prove that the Neumann data on side $C'$ is still
bounded independent of $\epsilon$, since our $\O( \epsilon)$ error
estimates are in terms of a priori Neumann data estimates on side $C'$.

\subsection{The upper bound of $\int_{ C'}|h\partial_{\nu}u|^2 dS$}
One of the assumptions we will use is that $\int_{{C'}}|h\partial_{\nu}u|^2 dS$ is bounded by a number that
is independent of $\epsilon$.
\begin{lemma}
  \label{L:1}
  Let $D \subset \reals^2$ be the domain from Theorem \ref{T:1}.  Then
  for $\epsilon>0$ sufficiently small, 
  there exists a number $\Gamma$, independent of $h$ and $\epsilon$ so
  that 
\[
\int_{C'} | h \partial_\nu u |^2 dS \leq \Gamma.
\]
\end{lemma}
The proof of Lemma \ref{L:1} is in Section \ref{SS:lemma} after  the
proofs of Theorems \ref{T:1} and \ref{T:2}.

\subsection{Proof of Theorem \ref{T:1} for acute triangles}
As shown in  Figure \ref{exampleplot}, sides $A, B, C'$ are listed in clockwise orientation. We use rectangular coordinates $(x,y)$ and orient our triangle such that the corner between sides $B$ and $C'$ is at the origin $(0,0)$, and the side $A$ is parallel to the $y$ axis.

Let $l$ be the segment on the $x$ axis that begins at $(0,0)$ and is
perpendicular to the side $A$. Write $A = A_1 \cup A_2$, where $A_1$
as the part of $A$ below the $x$ axis, and $A_2$ as the part above the
$x$ axis. Assume $a_1$ and $a_2$ be their lengths. We modify the acute
triangle by replacing its original side $C$ as described in the theorem. 

Specifically, we can parametrize side $B$ and $C'$ with respect to $x$:
\begin{equation*}
B=\{(x,y)\in \reals^{2}: y=-\frac{a_1}{l}x, 0\leq x\leq l\}
\end{equation*}
and 
\begin{equation*}
C'=\{(x,y)\in \reals^{2}:f(x)=\frac{a_2}{l}x +g(x),0\leq x\leq l\}
\end{equation*}
where $g(x)$ is restricted by the  two conditions from the theorem so that the function $f(x)$ is close to the original side of the acute triangle:
\begin{equation*}
|g(x)|\leq \epsilon 
\end{equation*}
\begin{equation*}
|g' (x)|\leq \epsilon.
\end{equation*}

Then the arclength parameters are $\gamma_A = 1$, 
\begin{equation*}
\gamma_B = \left(1 + (\frac{a_1}{l})^2\right)^\frac{1}{2} = \frac{(l^2 +a_1 ^2)^\frac{1}{2}}{l} = \frac{b}{l}
\end{equation*}
and
\begin{equation*}
\gamma_{C'} =(1+(f' (x))^2)^\frac{1}{2} 
\end{equation*}

We can then derive the unit tangent vectors:
\[
\tau_A = (0,1),
\]
\begin{equation*}
\tau_B = \left(1, -\frac{a_1}{l}\right)\gamma_B^{-1} = \left(\frac{l}{b},-\frac{a_1}{b}\right)
\end{equation*}
and 
\begin{equation*}
\tau_{C'} = (1,f\rq (x))\gamma_{C'}^{-1} = \left(\frac{1}{\sqrt{1+(f\rq (x))^2}},\frac{f\rq (x)}{\sqrt{1+(f\rq (x))^2}}\right).
\end{equation*}

From the unit tangent vectors, we find the outward unit normal vectors
to be:
\[
\nu_A = (1, 0),
\]
\begin{equation*}
\nu_B = \left(-\frac{a_1}{b},-\frac{l}{b}\right)
\end{equation*}
and 
\begin{equation*}
\nu_{C'} = \left(-\frac{f\rq (x)}{\sqrt{1+(f\rq (x))^2}},\frac{1}{\sqrt{1+(f\rq (x))^2}}\right)
\end{equation*}

The Dirichlet boundary conditions imply that the tangential derivatives of $u$ vanish on the boundary of the domain. That is
\begin{equation*}
\partial_y u = 0
\end{equation*}
on $A$, 
\begin{equation*}
\tau_B \cdot \nabla u = \frac{l}{b}\partial_x u-\frac{a_1}{b}\partial_y u = 0
\end{equation*}
on $B$, and
\begin{equation*}
\tau_{C'} \cdot \nabla u = \frac{1}{\sqrt{1+(f\rq (x))^2}}\partial_x u+\frac{f\rq (x)}{\sqrt{1+(f\rq (x))^2}}\partial_y u = 0
\end{equation*}
on ${C'}$. Rearranging, we get
\begin{equation*}
h\partial_x u = \frac{a_1}{l}h\partial_y u
\end{equation*}
on $B$, and
\begin{equation*}
h\partial_x u = -f\rq (x)h\partial_y u
\end{equation*}
on ${C'}$.

Next, we can relate $\partial_x$ and $\partial_y$ to $\partial_\nu$ on each side. Along $B$, we have
\begin{equation*} 
\begin{split}
h\partial_{\nu_B}u & = \nu_B\cdot h\nabla u \\
                  & = -\frac{a_1}{b}h\partial_x u -\frac{l}{b}h\partial_y u\\
& = \left(-\frac{a_1^2}{bl}-\frac{l}{b}\right)h\partial_y u \\
& = -\frac{b}{l}h\partial_y u
\end{split}
\end{equation*}
and thus
\begin{align*}
h\partial_y u & = -\frac{l}{b}h\partial_{\nu_B} u \\
h\partial_x u & = -\frac{a_1}{b}h\partial_{\nu_B} u
\end{align*}

Similarly, along ${C'}$, we have
\begin{equation*}
\begin{split}
h\partial_{\nu_{C'}} u & = \nu_{C'}\cdot h\nabla u \\
& = -\frac{f\rq (x)}{\sqrt{1+(f\rq (x))^2}}h\partial_x u + \frac{1}{\sqrt{1+(f\rq (x))^2}}h\partial_y u \\
& = (\frac{(f\rq (x))^2 +1}{\sqrt{1+(f\rq (x))^2}})h\partial_y u \\
\end{split}
\end{equation*}

Hence along ${C'}$ we have
\begin{align*}
h\partial_y u & = \frac{1}{\sqrt{1+(f\rq (x))^2}} h\partial_{\nu_{C'}} u \\
h\partial_x u & = -\frac{f\rq (x)}{\sqrt{1+(f\rq (x))^2}} h\partial_{\nu_{C'}} u
\end{align*}

Now consider the operator
\begin{equation*}
X = (x+m)\partial_x + (y+n)\partial_y
\end{equation*}
where $m, n$ are parameters that are independent of $x$ and $y$. The usual computation yields

\begin{equation*}
[-h^2 \Laplace - 1, X] = -2h^2 \Laplace
\end{equation*}

Then using the eigenfunction equation we have
\begin{equation*}
\begin{split}
\int_{D}([-h^2 \Laplace - 1, X]u)\overline{u}dV & = -2\int_{D}(h^2 \Laplace u)\overline{u}dV \\
&=\int_{D}2|u|^2 dV \\
&=2
\end{split}
\end{equation*}
since $u$ is normalized to the length of one.

Another way of calculation is
\begin{equation*}
\begin{split}
&\int_{D}([-h^2 \Laplace - 1, X]u)\overline{u}dV \\
& = \int_{D} ((-h^2 \Laplace - 1)Xu)\overline{u}dV - \int_{D}(X(-h^2\Laplace - 1)u)\overline{u}dV \\
&= \int_{D}((-h^2 \Laplace - 1)Xu)\overline{u}dV,
\end{split}
\end{equation*}
where we have used the eigenfunction equation in the last line.

Integrating by parts and applying the Green's theorem, we get
\begin{equation*}
\begin{split}
\int_{D}((-h^2 \Laplace - 1)Xu)\overline{u}dV &= \int_{D}(Xu)(-h^2 \Laplace - 1)\overline{u}dV \\
&- \int_{\partial D}(h\partial_{\nu}hXu)\overline{u}dS + \int_{\partial D}(hXu)(h\partial_{\nu} \overline{u})dS \\
&= \int_{\partial D}(hXu)(h\partial_{\nu} \overline{u})dS,
\end{split}
\end{equation*}
where we have used the Dirichlet boundary conditions in the last line.

Combining the results we have
\begin{equation*}
2 = \int_{\partial D}(hXu)(h\partial_{\nu} \overline{u})dS
\end{equation*}
which we can integrate on three sides separately. 

To simplify the notation, we define 
\begin{equation*}
I_A = \int_{A}|h\partial_{\nu}u|^2 dS
\end{equation*}
and similarly for $B$ and ${C'}$. 

Along $A$, we have
\begin{equation*}
\begin{split}
&\int_{A}(hXu)(h\partial_{\nu_A} \overline{u})dS \\
&= \int_{A}(((x+m)h\partial_x + (y+n)h\partial_y)u)(h\partial_{\nu_A} \overline{u})dS \\
&= (l+m)I_A
\end{split}
\end{equation*}
where $x = l$ on the side $A$.

On the side $B$, we substitute $y=-\frac{a_1}{l}x$ and get
\begin{equation*}
\begin{split}
&\int_{B}(hXu)(h\partial_{\nu_B} \overline{u})dS \\
&=\int_{B}\left(\left((x+m)h\partial_x + (-\frac{a_1}{l}x+n)h\partial_y\right)u\right)(h\partial_{\nu_B} \overline{u})dS \\
&=\int_{B}\left(\left(\left((x+m)\left(-\frac{a_1}{b}\right) + \left(-\frac{a_1}{l}x+n\right)\left(-\frac{l}{b}\right)\right)h\partial_{\nu_B}\right)u\right)(h\partial_{\nu_B} \overline{u})dS \\
&=\int_{B}\left(\left(-\frac{a_1}{b}m-\frac{l}{b}n\right)h\partial_{\nu_B}u\right)(h\partial_{\nu_B} \overline{u})dS \\
&=\left(-\frac{a_1}{b}m-\frac{l}{b}n\right)I_B
\end{split}
\end{equation*}

On ${C'}$, we substitute $y=f(x)$ and get
\begin{equation*}
\begin{split}
&\int_{{C'}}(hXu)(h\partial_{\nu_{C'}} \overline{u})dS \\
&=\int_{{C'}}(((x+m)h\partial_x + (y+n)h\partial_y)u)(h\partial_{\nu_{C'}} \overline{u})dS \\
&=\int_{{C'}}\left(\left(-(x+m)\frac{f\rq (x)}{\sqrt{1+f\rq (x)^2}}+(f(x)+n)\frac{1}{\sqrt{1+f\rq (x)^2}}\right)\partial_{\nu_{C'}} u\right)(h\partial_{\nu_{C'}} \overline{u})dS \\
\end{split}
\end{equation*}

Hence, summing up the integrations along the three sides, we have
\begin{equation*}
\begin{split}
&\int_{\partial D}(hXu)(h\partial_{\nu_A} \overline{u})dS \\
&= (l+m)I_A+\left(-\frac{a_1}{b}m-\frac{l}{b}n\right)I_B \\
&+ \int_{{C'}}\left(\left(-(x+m)\frac{f\rq (x)}{\sqrt{1+f\rq (x)^2}}+(f(x)+n)\frac{1}{\sqrt{1+f\rq (x)^2}}\right)h\partial_{\nu_{C'}} u\right)(h\partial_{\nu_{C'}} \overline{u})dS \\
&=2
\end{split}
\end{equation*}

First, as $m$ and $n$ are independent parameters, we can set $m=n=0$, which yields
\begin{equation}
  \label{E:IA}
lI_A + \int_{{C'}}\left(\left(\frac{-xf\rq (x) +f(x) }{\sqrt{1+f\rq (x)^2}}\right)h\partial_{\nu_{C'}} u\right)(h\partial_{\nu_{C'}} \overline{u})dS =2
\end{equation}

Additionally, we can differentiate with respect to $m$:
\begin{equation}
  \label{E:IB}
I_A - \frac{a_1}{b}I_B - \int_{{C'}} \left(\frac{f\rq (x)}{\sqrt{1+f\rq (x)^2}}h\partial_{\nu_{C'}} u\right)(h\partial_{\nu_{C'}} \overline{u})dS = 0
\end{equation}
and with respect to $n$:
\begin{equation}
  \label{E:IC}
- \frac{l}{b}I_B + \int_{{C'}} \left(\left(\frac{1}{\sqrt{1+f\rq (x)^2}}\right)h\partial_{\nu_{C'}} u\right)(h\partial_{\nu_{C'}} \overline{u})dS = 0.
\end{equation}

In equation \eqref{E:IA}, observe that
\begin{equation*}
\begin{split}
|-xf\rq (x) +f(x)| &= \left|-x\left(\frac{a_2}{l}+g\rq (x)\right) + \frac{a_2}{l}x+g(x)\right| \\
&= |g(x) - xg\rq (x)| \\
&\leq |g(x)| + |xg\rq (x)| \\
&\leq \epsilon + l\epsilon \\
& = (l+1)\epsilon
\end{split}
\end{equation*}
as $x\leq l$. Plugging this in equation \eqref{E:IA} yields
\begin{equation*}
\begin{split}
2 &= lI_A + \int_{{C'}}\left(\left(\frac{-xf\rq (x) +f(x) }{\sqrt{1+f\rq (x)^2}}\right)h\partial_{\nu_{C'}} u\right)(h\partial_{\nu_{C'}} \overline{u})dS \\
&\leq lI_A + \int_{{C'}}\left(\left(\frac{|-xf\rq (x) +f(x)|}{\sqrt{1+f\rq (x)^2}}\right)h\partial_{\nu_{C'}} u\right)(h\partial_{\nu_{C'}} \overline{u})dS \\
&\leq lI_A + \int_{{C'}}\left(\left(\frac{(l+1)\epsilon}{\sqrt{1+f\rq (x)^2}}\right)h\partial_{\nu_{C'}} u\right)(h\partial_{\nu_{C'}} \overline{u})dS \\
&=lI_A + \int_{{C'}}\frac{(l+1)\epsilon}{\sqrt{1+f\rq (x)^2}}|h\partial_{\nu}|^2 dS \\
&\leq lI_A + \int_{{C'}}(l+1)\epsilon|h\partial_{\nu}|^2 dS \\
& = lI_A + \beta \epsilon
\end{split}
\end{equation*}
where $\beta =\int_{{C'}}(l+1)|h\partial_{\nu}|^2 dS = (l+1)I_{C'}$. In the next section, we will show that $\beta$ is finite and bounded by a number that is independent of $\epsilon$. 

Thus, we find the lower bound of $I_A$ to be
\begin{equation*}
\frac{2}{l}-\frac{\beta \epsilon}{l} \end{equation*}

The upper bound of $I_A$ can be found in a similar way to get
\begin{equation*}
I_A \leq  \frac{2}{l}+\frac{\beta \epsilon}{l}.
\end{equation*}

Hence, we find the range of $I_A$ to be
\begin{equation*}
\frac{2}{l}-\frac{\beta \epsilon}{l} \leq I_A \leq  \frac{2}{l}+\frac{\beta \epsilon}{l}.
\end{equation*}
Now, comparing to the original triangle $T$, we have 
\begin{equation*}
\frac{2}{l} = \frac{a}{al/2} =\frac{a}{\text{Area}(T)},
\end{equation*}
and the perturbation $g$ changes the area by a factor controlled by
$\epsilon$:
\[
\text{Area}(D) = \text{Area}(T) + \O(\epsilon).\]
Hence
\[
\frac{2}{l} = \frac{a}{\text{Area}(D)} + \O(\epsilon).
\]
In other words, we find
\begin{equation*}
I_A = \frac{a}{Area(D)} + \O(\epsilon).
\end{equation*}
Note that when $\epsilon = 0$, the domain D would be a triangle, and
$I_A = \frac{a}{Area(D)}$, so this is consistent with the results in
\cite{Chr-tri}.

Next, we substitute $f\rq (x) = \frac{a_2}{l}+g\rq (x) $ in equation
\eqref{E:IB} and use equation \eqref{E:IC} to get 
\begin{equation*}
\begin{split}
&I_A - \frac{a_1}{b}I_B - \int_{{C'}} \left(\frac{f\rq (x)}{\sqrt{1+(f\rq (x))^2}}h\partial_{\nu_{C'}} u\right)(h\partial_{\nu_{C'}} \overline{u})dS \\
&= I_A - \frac{a_1}{b}I_B -\int_{{C'}} \left(\frac{\frac{a_2}{l}+g\rq (x)}{\sqrt{1+(f\rq (x))^2}}h\partial_{\nu_{C'}} u\right)(h\partial_{\nu_{C'}} \overline{u})dS \\
&= I_A - \frac{a_1}{b}I_B -\frac{l}{b}\frac{a_2}{l}I_B - \int_{{C'}} \left(\frac{g\rq (x)}{\sqrt{1+(f\rq (x))^2}}h\partial_{\nu_{C'}} u\right)(h\partial_{\nu_{C'}}\overline{u})dS \\
&= I_A - \frac{a}{b}I_B - \int_{{C'}} \left(\frac{g\rq (x)}{\sqrt{1+f\rq (x)^2}}h\partial_{\nu_{C'}} u\right)(h\partial_{\nu_{C'}}\overline{u})dS \\
& = 0
\end{split}
\end{equation*}
and thus
\begin{equation*}
I_B = \frac{b}{a}I_A - \frac{b}{a}\int_{{C'}} \left(\frac{g\rq (x)}{\sqrt{1+f\rq (x)^2}}h\partial_{\nu_{C'}} u\right)(h\partial_{\nu_{C'}}\overline{u})dS
\end{equation*}

As we assume that $|g\rq (x)|$ is bounded by $\epsilon$, we can find
the upper and lower bound of $I_B$ as we did for $I_A$:
\begin{equation*}
\begin{split}
I_B &=\frac{b}{a}I_A - \frac{b}{a}\int_{{C'}} \left(\frac{g\rq
  (x)}{\sqrt{1+(f\rq (x))^2}}h\partial_{\nu_{C'}}
u\right)(h\partial_{\nu_{C'}}\overline{u})dS \\
& = \frac{b}{a}I_A + \O(\epsilon) \\
& =\left( \frac{b}{a} \right) \frac{a}{\text{Area}(D)} + \O(\epsilon) \\
& = \frac{b}{\text{Area}(D)} + \O(\epsilon).
\end{split}
\end{equation*}




Finally, we plug in the range of $I_B$ to equation \eqref{E:IC} and find
\begin{align}
  \int_{{C'}}  \left(\frac{1}{\sqrt{1+(f\rq
      (x))^2}}h\partial_{\nu_{C'}} u\right)(h\partial_{\nu_{C'}}
  \overline{u})dS & = 
  \frac{l}{b} I_B  \\
  & = \frac{l}{\text{Area}(D)} + \O(\epsilon).
  \label{E:IC2}
\end{align}

In order to find the range of $\frac{1}{\sqrt{1+f\rq (x)^2}}$, we use
our assumption that $|g\rq (x)|$ is bounded by $\epsilon$ and get
\begin{equation*}
\begin{split}
(f\rq (x))^2 &= \left(\frac{a_2}{l}+g\rq (x)\right)^2\\
&= \left(\frac{a_2}{l}\right)^2 + \frac{2a_2}{l}g\rq (x) + g\rq (x)^2\\
&= \left(\frac{a_2}{l}\right)^2 + \alpha 
\end{split}
\end{equation*}
where $\alpha$ is a function of $g\rq (x)$, $\alpha = \O(\epsilon)$
for $\epsilon$ small. 

Next, we have
\begin{equation*}
\begin{split}
\frac{1}{\sqrt{1+f\rq (x)^2}} &=\frac{1}{\sqrt{1+(\frac{a_2}{l})^2 +\alpha }}\\
&= \frac{1}{\sqrt{(1+(\frac{a_2}{l})^2)}}\frac{1}{\sqrt{(1+\frac{\alpha}{1+(\frac{a_2}{l})^2})}}\\
&=\frac{l}{c}\frac{1}{\sqrt{1+\frac{l}{c}\alpha}} \\
& = \frac{l}{c} + \O(\epsilon),
\end{split}
\end{equation*}
where $c = (a_2^2 + l^2)^{1/2}$ is the length of side $C$ before
deforming it to side $C'$.



Plugging this result back to equation \eqref{E:IC2} yields
\begin{equation}
  \begin{split}
\frac{2}{a} + \O(\epsilon) & = 
\int_{{C'}}\left(\frac{1}{\sqrt{1+f\rq (x)^2}}h\partial_{\nu_{C'}}
u\right)(h\partial_{\nu_{C'}}\overline{u})dS \\
&= \frac{l}{c}(1+ \O(\epsilon))\int_{{C'}}|h\partial_{\nu_{C'}}u|^2dS. 
  \end{split}
  \label{E:IC3}
\end{equation}

We now use again that
\[
\text{Area}(D) = \text{Area}(T) + \O(\epsilon) = \frac{al}{2} +
\O(\epsilon)
\]
where $T$ is the original triangle, and that the length of $C'$ is
$l(C') = c + \O(\epsilon)$.  
Therefore, we have from \eqref{E:IC3}
\begin{align*}
\int_{{C'}}|h\partial_{\nu_{C'}}u|^2dS & = \frac{c}{l}
 \frac{2}{a} + \O(\epsilon)\\
 & = \frac{c}{al/2} + \O(\epsilon) \\
 & = \frac{l({C'})}{\text{Area}(D)} + \O(\epsilon).
 \end{align*}
This proves the theorem in the case of an acute or right triangle.
  


\subsection{Proof of Theorem \ref{T:1} for obtuse triangles}
The proof of obtuse triangles is nearly the same, with a few changes in the signs. The set up for obtuse triangles is shown in Figure \ref{exampleplot2}.

We can parametrize $B$ and $C'$ with respect to $x$:
\begin{equation*}
B=\{(x,y)\in \reals^{2}:f(x)=\frac{a_1}{l}x,0\leq x\leq l\}
\end{equation*}
and 
\begin{equation*}
C'=\{(x,y)\in \reals^{2}:f(x)=\frac{a_2
+a_1}{l}x + g(x),0\leq x\leq l\}
\end{equation*}

Doing similar computations as before 
we find
\begin{align*}
h\partial_y u & = -\frac{l}{b}h\partial_{\nu_B} u \\
h\partial_x u & = \frac{a_1}{b}h\partial_{\nu_B} u
\end{align*}
along $B$, and 
\begin{align*}
h\partial_y u & = \frac{1}{\sqrt{1+f\rq (x)^2}} h\partial_{\nu_{C'}} u \\
h\partial_x u & = -\frac{f\rq (x)}{\sqrt{1+f\rq (x)^2}} h\partial_{\nu_{C'}} u
\end{align*}
along ${C'}$.

Following the commutator computation as in the acute case, and plug in the equation of side $B$ and ${C'}$, we have
\begin{equation*}
\begin{split}
2 &= (l+m)I_A+\left(\frac{a_1}{b}m-\frac{l}{b}n\right)I_B \\ 
&+ \int_{{C'}}\left(\left(-(x+m)\frac{f\rq (x)}{\sqrt{1+f\rq (x)^2}}+(f(x)+n)\frac{1}{\sqrt{1+f\rq (x)^2}}\right)h\partial_{\nu_{C'}} u\right)(h\partial_{\nu_{C'}} \overline{u})dS 
\end{split}
\end{equation*}
Differentiating with respect to $m$ and $n$ we get
\begin{equation}
I_A + \frac{a_1}{b}I_B - \int_{{C'}} \left(\frac{f\rq (x)}{\sqrt{1+f\rq
    (x)^2}}h\partial_{\nu_{C'}} u\right)(h\partial_{\nu_{C'}}
\overline{u})dS = 0
\label{E:IAo}
\end{equation}
and 
\begin{equation}
- \frac{l}{b}I_B + \int_{{C'}} \left(\left(\frac{1}{\sqrt{1+f\rq
    (x)^2}}\right)h\partial_{\nu_{C'}} u\right)(h\partial_{\nu_{C'}}
\overline{u})dS = 0.
\label{E:ICo}
\end{equation}

Again, if we set $m = n = 0$, we find the range of $I_A$ to be
\begin{equation*}
  I_A =  \frac{2}{l}+ \O(\epsilon) = \frac{a_2}{\text{Area}(D)} + \O(\epsilon).
\end{equation*}

Next, plugging in $f'(x)$ in equation \eqref{E:IAo} and using equation \eqref{E:ICo}, we have
\begin{equation*}
\begin{split}
&I_A + \frac{a_1}{b}I_B - \int_{{C'}} \left(\frac{f\rq (x)}{\sqrt{1+f\rq (x)^2}}h\partial_{\nu_{C'}} u\right)(h\partial_{\nu_{C'}} \overline{u})dS \\
&= I_A + \frac{a_1}{b}I_B -\int_{{C'}} \left(\frac{\frac{a_1 + a_2}{l}+g\rq (x)}{\sqrt{1+f\rq (x)^2}}h\partial_{\nu_{C'}} u\right)(h\partial_{\nu_{C'}} \overline{u})dS \\
&= I_A + \frac{a_1}{b}I_B -\frac{l}{b}\frac{a_1 + a_2}{l}I_B - \int_{{C'}} \left(\frac{g\rq (x)}{\sqrt{1+f\rq (x)^2}}h\partial_{\nu_{C'}} u\right)(h\partial_{\nu_{C'}}\overline{u})dS \\
&= I_A - \frac{a_2}{b}I_B - \int_{{C'}} (\frac{g\rq (x)}{\sqrt{1+f\rq (x)^2}}h\partial_{\nu_{C'}} u)(h\partial_{\nu_{C'}}\overline{u})dS \\
& = 0
\end{split}
\end{equation*}
which is the same equation we have for acute triangles. 

Therefore, using the range of $I_A$ and the same estimates as in the
acute case, we find the range of $I_B$ to be
\begin{equation*}
I_B = \frac{2b}{a_2l} + \O(\epsilon) = \frac{b}{\text{Area}(D)} + \O(\epsilon).
\end{equation*}

Finally, using equation \eqref{E:ICo} and following the computation above, the range of $I_{C'}$ is the same as that of acute triangles:
\begin{align*}
  I_{C'} & =\frac{2c}{a_2l} + \O(\epsilon) = \frac{l(C')}{ \text{Area}(D)} + \O(\epsilon).
\end{align*}

\section{Proof of Theorem \ref{T:2}}
We now proceed with the proof of Theorem \ref{T:2}.  It naturally is
very similar to that of Theorem \ref{T:1} so we just point out some of
the main differences.  
\begin{proof}
  With the same vector field $X = (x + m) \p_x + (y + n ) \p_y$, the calculation of the commutator alone tells us that
\begin{equation*}
[-h^2 \Laplace - 1, X] = -2h^2 \Laplace = 2(-h^2 \Laplace +w(x,y)) -2w(x,y)
\end{equation*}
and
\begin{equation*}
\begin{split}
\int_{T}([-h^2 \Laplace - 1, X]u)\overline{u}dV & = 2\int_{T}((-h^2
\Laplace +w(x, y)) u)\overline{u}dV - 2\int_{T}(w(x, y)u)\overline{u}dV \\
&=2- 2\int_{T}(w(x,y)u)\overline{u}dV 
\end{split}
\end{equation*}

Because $|w(x,y)|<\epsilon$, we have
\begin{equation*}
\begin{split}
\left|\int_{T}([-h^2 \Laplace - 1, X]u)\overline{u}dV\right|  
&=\left|2- 2\int_{T}(w(x,y)u)\overline{u}dV\right|\\
&\geq 2 - 2\epsilon\int_{T}|u|^2 dV\\
&= 2 - 2\epsilon
\end{split}
\end{equation*}
and 
\begin{equation*}
\left|\int_{T}([-h^2 \Laplace - 1, X]u)\overline{u}dV\right|\leq 2 + 2\epsilon.
\end{equation*}

On the other hand, we have
\begin{equation}
\begin{split}
&\int_{T}([-h^2 \Laplace - 1, X]u)\overline{u}dV \\
& = \int_{T} ((-h^2 \Laplace - 1)Xu)\overline{u}dV - \int_{T}(X(-h^2\Laplace - 1)u)\overline{u}dV \\
&= \int_{T}((-h^2 \Laplace - 1)Xu)\overline{u}dV + \int_{T}(X(w(x,y)u)\overline{u}dV
\end{split}
\label{E:IBP-w-1}
\end{equation}

Integrating by parts, we have
\begin{equation}
\begin{split}
\int_{T}((-h^2 \Laplace - 1)Xu)\overline{u}dV &= \int_{T}(Xu)(-h^2 \Laplace - 1)\overline{u}dV \\
&- \int_{\partial T}(h\partial_{\nu}hXu)\overline{u}dS + \int_{\partial T}(hXu)(h\partial_{\nu} \overline{u})dS \\
&= -\int_{T}(Xu)w(x, y)\overline{u}dV + \int_{\partial T}(hXu)(h\partial_{\nu} \overline{u})dS.
\end{split}
\label{E:IBP-w-2}
\end{equation}
The last term in \eqref{E:IBP-w-1} is computed:
\begin{equation*}
\begin{split}
\int_{T}(X(w(x,y)u)\overline{u}dV &= \int_{T}((Xw(x,y))u) + (w(x,y)Xu))\overline{u} dV\\
&= \int_{T}((Xw(x,y))u)\overline{u}dV + \int_{T}(w(x,y)Xu)\overline{u}dV.
\end{split}
\end{equation*}

Combining this with \eqref{E:IBP-w-1} and \eqref{E:IBP-w-2}, we have
\begin{equation*}
\begin{split}
&\int_{T}([-h^2 \Laplace - 1, X]u)\overline{u}dV \\
&= \int_{T}((-h^2 \Laplace - 1)Xu)\overline{u}dV + \int_{T}(X(w(x,y)u)\overline{u}dV \\
&= \int_{\partial T}(hXu)(h\partial_{\nu} \overline{u})dS + \int_{T}((Xw(x,y))u)\overline{u}dV.
\end{split}
\end{equation*}
With the condition that $|\nabla w(x,y|) \leq \epsilon$, we get
\begin{equation*}
\begin{split}
&\left|\int_{T}((Xw(x,y))u)\overline{u}dV\right| \\
&= \left|\int_{T} (((x+m)\partial_x w(x,y) + (y + n) \p_y w)x,y))u)\overline{u} dV\right|\\
  &= \O(\epsilon) + m \O(\epsilon) + n \O(\epsilon).
\end{split}
\end{equation*}

Combining the results together, we have
\begin{equation*}
\begin{split}
2+ \O(\epsilon) + m \O(\epsilon) + n \O(\epsilon) &= \int_{T}([-h^2 \Laplace - 1, X]u)\overline{u}dV\\
& = \int_{\partial T}(hXu)(h\partial_{\nu}\overline{u})dS. 
\end{split}
\end{equation*}
The rest of the proof proceeds exactly as the proof of Theorem \ref{T:1}.
\end{proof}

\section{Proof of Lemma \ref{L:1}}
\label{SS:lemma}

\begin{proof}
To prove the Lemma, first we consider the vector field
\begin{equation*}
X = y\partial_y
\end{equation*}
and the usual computation yields
\begin{equation*}
[-h^2 \Laplace-1, y\partial_y] = -2h^2\partial_y^2
\end{equation*}

Then the integration yields
\begin{equation*}
\int_D ([-h^2 \Laplace-1, y\partial_y]u)\overline{u} dV = -2\int_D (h^2\partial_y^2u)\overline{u} dV
\end{equation*}

Moreover, observe that
\begin{equation*}
\begin{split}
-2\int_D (h^2\partial_y^2u)\overline{u} dV &= 2\int_{D}|h\partial_y u|^2 dV \\
&\leq 2\int_{D}(|h\partial_y u|^2+|h\partial_x u|^2) dV \\
&= 2\int_{D} |u|^2 dV \\
&= 2
\end{split}
\end{equation*}

On the other hand, if we integrate by parts and using the boundary conditions, we have
\begin{equation*}
\begin{split}
\int_D ([-h^2 \Laplace-1, y\partial_y]u)\overline{u} dV &= \int_D ((-h^2 \Laplace-1)y\partial_yu)\overline{u} dV \\
&=\int_{\partial D}(yh\partial_yu)(h\partial_{\nu}\overline{u})dS
\end{split}
\end{equation*}

Hence, together we have
\begin{equation*}
-2\int_D (h^2\partial_y^2u)\overline{u} dV = \int_{\partial D}(yh\partial_yu)(h\partial_{\nu}\overline{u})dS \leq 2
\end{equation*}

Since $A$ is vertical, we have 
\begin{equation*}
\partial_y u =0
\end{equation*}
on $A$,
\[
\partial_y u = - \frac{l}{b} \partial_\nu
\]
on $B$ (in both the acute and obtuse cases), and
\begin{equation*}
\begin{split}
\partial_y u &= \frac{1}{\sqrt{1+(f'(x))^2}}\partial_{\nu} u \\
&= \gamma^{-1} \partial_{\nu} u
\end{split}
\end{equation*}
 on $C'$, where $\gamma = {\sqrt{1+(f'(x))^2}}$ is the arclength element.   Substituting $\partial_y u$, we have
\begin{equation}
\begin{split}
&\int_{\partial D}(yh\partial_yu)(h\partial_{\nu}\overline{u})dS \\
&= \int_{C'} f(x)\gamma^{-1} |h\partial_{\nu} u|^2 dS + \int_{B}
    \left(-\frac{l}{b}\right)\left(\mp \frac{a_2}{l} x \right) |h\partial_{\nu} u|^2 dS \\
    &= \int_{C'} f(x)\gamma^{-1} |h\partial_{\nu} u|^2 dS \pm
    \frac{a_2}{b} \int_B x | h \partial_\nu u |^2 dS
    \\ 
&\leq 2,
\end{split}
\label{E:a-priori-1}
\end{equation}
where the $\pm$ sign corresponds to the acute/obtuse cases.

While this is close to what we intend to prove, we should be careful
because $f(x)$ approaches zero as $x$ goes to zero.  There is also a
potential problem in the obtuse case because of the sign change on the
$B$ integral.

\noindent {\bf Acute case:}  
Observe that in \eqref{E:a-priori-1} the function $a_2x/b \geq 0$, so
in the acute case, we have
\begin{align*}
  & \left| \int_{C'} f(x)\gamma^{-1} |h\partial_{\nu} u|^2 dS \right| \\
  & \quad =  \int_{C'} f(x)\gamma^{-1} |h\partial_{\nu} u|^2 dS \\
  & \quad \leq  \int_{C'} f(x)\gamma^{-1} |h\partial_{\nu} u|^2 dS
  +
    \frac{a_2}{b} \int_B x | h \partial_\nu u |^2 dS
    \\ 
    & \quad \leq 2
\end{align*}
so to show that $\int_{C'} | h \p_\nu u |^2 dS$ is bounded we only
need to estimate the integral for $x$ near $0$.

Fix $\delta > 0$ such that $\epsilon < \frac{a\delta}{2b}$ and $\delta
\gg \epsilon$.    Using our restrictions on $f(x)$, we can find the lower bound of $f(x)$ when $x \geq \delta$:
\begin{equation*}
\begin{split}
f(x) &= \frac{a}{b}x + g(x) \\
&\geq\frac{a}{b}x - |g(x)| \\
&\geq\frac{a}{b}x -\epsilon\\
&\geq \frac{a\delta}{b} -\epsilon \\
&\geq \frac{a\delta}{2b}.
\end{split}
\end{equation*}
For $\delta \leq x \leq l$, we have  an upper bound for $\gamma^{-1}$:
\begin{equation*}
\begin{split}
\gamma &={\sqrt{1+(f'(x))^2}} \\
&= {\sqrt{1+\left(\frac{a}{b} + g'(x)\right)^2}} \\
&\leq {\sqrt{1+\left(\frac{a}{b} + \epsilon\right)^2}}\\
&= : \gamma_0.
\end{split}
\end{equation*}
We observe that then
\begin{align*}
\gamma^{-1} \geq \gamma_0^{-1} & = \frac{1}{\sqrt{1 + \frac{a^2}{b^2}}}
+ \O(\epsilon) \\
& \geq \frac{1}{2\sqrt{1 + \frac{a^2}{b^2}}}
\end{align*}
for $\epsilon>0$ sufficiently small.

Hence, substituting $f(x)$ and $\gamma^{-1}$ with their lower bounds, we have
\begin{equation}
\begin{split}
2 &\geq \int_{C'} f(x)\gamma^{-1} |h\partial_{\nu} u|^2 dS
  \\
& \geq \gamma_0^{-1} \frac{a\delta}{2b} \int_{C' \cap \{ x \geq \delta
  \}} | h \partial_\nu u |^2 dS  \\
& \geq \frac{a \delta}{4b \sqrt{1 + \frac{a^2}{b^2}}}\int_{C' \cap \{ x \geq \delta
  \}} | h \partial_\nu u |^2 dS.
\end{split}
\label{E:I-delta1}
\end{equation}
That means that \eqref{E:I-delta1} implies
\[
\int_{C' \cap \{ x \geq \delta
  \}} | h \partial_\nu u |^2 dS = \O_\delta (1)
\]
independent of $\epsilon$ provided $\epsilon$ is sufficiently small.

For $x < \delta$, consider another function $\psi(x)$, which has value
one on $x \leq \delta$ and monotonically decreases to zero for $x \geq
2\delta$.   Such a function $\psi$ is depicted in Figure \ref{F:phi}.

    \begin{figure}
\hfill
\centerline{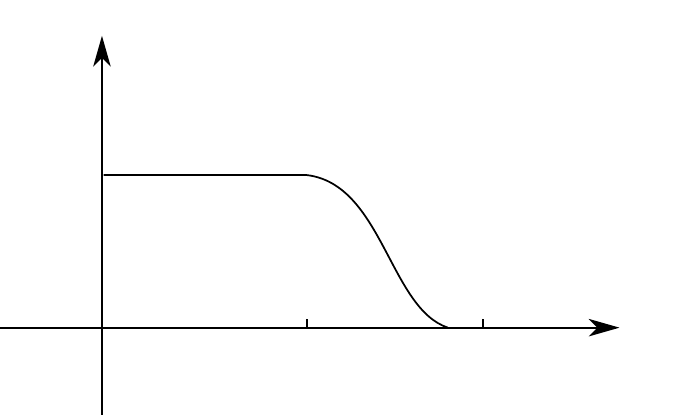}
\caption{\label{F:phi}  The function $\psi$}
\hfill
\end{figure}


Then computing the commutator we have
\begin{equation*}
[-h^2 \Laplace -1, \psi\partial_x]=-2\psi'h^2\partial_x^2-h\psi''h\partial_x
\end{equation*}
Thus we have
\begin{equation*}
\begin{split}
\left|\int_{D}([-h^2 \Laplace -1, \psi\partial_x]u)\overline{u}dV\right| &= \left|\int_{D}(2\psi'h^2\partial_x^2u)\overline{u}+(h\psi''h\partial_xu)\overline{u}dV\right| \\
&=\left|\int_{D}(-2h\partial_x u h\partial_x\psi' \overline{u})dV+\int_{D}(h\psi''h\partial_xu)\overline{u}dV\right| \\
&= \left|\int_{D}(-h\partial_x uh\psi''\overline{u}-2\psi'|h\partial_x u|^2)dV\right| \\
&\leq \sup(|2\psi'|,|\psi''|)\int_{D}(|h\partial_x
u||hu|+|h\partial_x u|^2)dV \\
& = \O_\delta (1).
\end{split}
\end{equation*}
Here the implicit constant in  the $\O_\delta(1)$ depends on our fixed
$\delta$, 
but not on $\epsilon \ll \delta$.

On the other hand, if we integrate by parts, we have
\begin{equation*}
\int_{D}([-h^2 \Laplace -1, \psi\partial_x]u)\overline{u}dV =
\int_{\partial D} \psi h\partial_xuh\partial_{\nu}\overline{u}dS = \O_\delta(1).
\end{equation*}

Since 
\begin{align*}
  \partial_x u & = - \frac{f'}{\gamma} \partial_\nu u 
\end{align*}
on $C'$ and we have already computed
\begin{align*}
 \left| \frac{f'}{\gamma} \right| & = \left| \frac{ \frac{a_2}{l} + g'
 }{\gamma} \right| \\
 & \geq \frac{a_2}{4l\sqrt{1 + \frac{a^2}{b^2}}}
\end{align*}
for $\epsilon$ sufficiently small, 
we have
\begin{equation*}
  \begin{split}
\int_{C' \cap \{0 \leq x \leq \delta\}} | h \p_\nu u |^2 dS & \leq
\int_{C'} \psi(x) | h \p_\nu u |^2 dS \\
& \leq \frac{4l\sqrt{1 + \frac{a^2}{b^2}}}{a_2}\int_{C'} \psi
\frac{f'}{\gamma} | h \p_\nu u |^2 dS \\
& = \O_\delta(1).
  \end{split}
  \end{equation*}


Combining with \eqref{E:I-delta1}, we have
\[
\int_{C'} | h \p_\nu u |^2 dS 
= \O(1),
\]
which proves the Lemma in this cases.

\noindent {\bf Case of obtuse triangle:}  In this case we have to be
slightly more careful.  Consider the vector field
\[
X = \left( y - \frac{a_2}{l} x\right) \p_y.
\]
We have $[-h^2 \Laplace -1 , X] = -2 h^2 \p_y^2 + 2 \frac{a_2}{l} h
\p_x h \p_y$.  The interior estimates are similar, so that
\[
\int_D ([-h^2 \Laplace -1 , X] u) \overline{u} dV = \O(1).
\]
Then the vector field $X$ vanishes when $y = \frac{a_2}{l} x$.
Further, since $X$ is tangential on side $A$, fixing a $\delta \gg
\epsilon$, the same argument as in
the acute case gives
\[
\int_{C'\cap \{ \delta \leq x \leq l \}} | h \p_\nu u |^2 dS = \O(1).
\]
For the set $\{ 0 \leq x \leq \delta\}$, we use the vector field $Y =
\psi(x)\left( \p_x + 
\frac{a_2}{l} \p_y \right)$.  Then $Y = 0$ on $A$ since $\psi = 0$
there, and $Y$ is tangential to $B$, so $Y u = 0$ on $B$.

On $C'$, we have
\[
h Y u = \psi(x) \left( - \frac{f'}{\gamma} + \frac{a_1}{l \gamma}
\right) h \p_\nu u.
\]
Since in the obtuse case we  have $f' = (a_1 + a_1)/l + g'$, 
\[
\frac{f' - \frac{a_1}{l}}{\gamma} = \frac{a_2 }{l \gamma} +
\O(\epsilon) \geq \frac{a_2 }{2 l \gamma}
\]
independent of $\epsilon$ sufficiently small.  Using the previously
established estimates on $\gamma$, the rest of the proof follows
exactly as in the acute case.

\end{proof}

\bibliographystyle{alpha}
\bibliography{HC-bib}

\end{document}

%% file: acute-1.pdf_tex
\begingroup%
  \makeatletter%
  \providecommand\color[2][]{%
    \errmessage{(Inkscape) Color is used for the text in Inkscape, but the package 'color.sty' is not loaded}%
    \renewcommand\color[2][]{}%
  }%
  \providecommand\transparent[1]{%
    \errmessage{(Inkscape) Transparency is used (non-zero) for the text in Inkscape, but the package 'transparent.sty' is not loaded}%
    \renewcommand\transparent[1]{}%
  }%
  \providecommand\rotatebox[2]{#2}%
  \newcommand*\fsize{\dimexpr\f@size pt\relax}%
  \newcommand*\lineheight[1]{\fontsize{\fsize}{#1\fsize}\selectfont}%
  \ifx\svgwidth\undefined%
    \setlength{\unitlength}{276.42552567bp}%
    \ifx\svgscale\undefined%
      \relax%
    \else%
      \setlength{\unitlength}{\unitlength * \real{\svgscale}}%
    \fi%
  \else%
    \setlength{\unitlength}{\svgwidth}%
  \fi%
  \global\let\svgwidth\undefined%
  \global\let\svgscale\undefined%
  \makeatother%
  \begin{picture}(1,0.52912037)%
    \lineheight{1}%
    \setlength\tabcolsep{0pt}%
    \put(0,0){\includegraphics[width=\unitlength,page=1]{acute-1.pdf}}%
    \put(0.94628638,0.26745636){\color[rgb]{0,0,0}\makebox(0,0)[lt]{\lineheight{1.25}\smash{\begin{tabular}[t]{l}$x$\end{tabular}}}}%
    \put(0.78071802,0.35788212){\color[rgb]{0,0,0}\makebox(0,0)[lt]{\lineheight{1.25}\smash{\begin{tabular}[t]{l}$A$\end{tabular}}}}%
    \put(0.68265069,0.32858924){\color[rgb]{0,0,0}\makebox(0,0)[lt]{\lineheight{1.25}\smash{\begin{tabular}[t]{l}$a_2$\end{tabular}}}}%
    \put(0.68519788,0.17957786){\color[rgb]{0,0,0}\makebox(0,0)[lt]{\lineheight{1.25}\smash{\begin{tabular}[t]{l}$a_1$\end{tabular}}}}%
    \put(0.48014794,0.1209921){\color[rgb]{0,0,0}\makebox(0,0)[lt]{\lineheight{1.25}\smash{\begin{tabular}[t]{l}$B$\end{tabular}}}}%
    \put(0.44448718,0.42410945){\color[rgb]{0,0,0}\makebox(0,0)[lt]{\lineheight{1.25}\smash{\begin{tabular}[t]{l}$C'$\end{tabular}}}}%
    \put(0.18339871,0.50434643){\color[rgb]{0,0,0}\makebox(0,0)[lt]{\lineheight{1.25}\smash{\begin{tabular}[t]{l}$y$\end{tabular}}}}%
    \put(0.28401327,0.37189177){\color[rgb]{0,0,0}\makebox(0,0)[lt]{\lineheight{1.25}\smash{\begin{tabular}[t]{l}$y=f(x)$\end{tabular}}}}%
  \end{picture}%
\endgroup%

%% file: obtuse-1.pdf_tex
\begingroup%
  \makeatletter%
  \providecommand\color[2][]{%
    \errmessage{(Inkscape) Color is used for the text in Inkscape, but the package 'color.sty' is not loaded}%
    \renewcommand\color[2][]{}%
  }%
  \providecommand\transparent[1]{%
    \errmessage{(Inkscape) Transparency is used (non-zero) for the text in Inkscape, but the package 'transparent.sty' is not loaded}%
    \renewcommand\transparent[1]{}%
  }%
  \providecommand\rotatebox[2]{#2}%
  \newcommand*\fsize{\dimexpr\f@size pt\relax}%
  \newcommand*\lineheight[1]{\fontsize{\fsize}{#1\fsize}\selectfont}%
  \ifx\svgwidth\undefined%
    \setlength{\unitlength}{276.42552567bp}%
    \ifx\svgscale\undefined%
      \relax%
    \else%
      \setlength{\unitlength}{\unitlength * \real{\svgscale}}%
    \fi%
  \else%
    \setlength{\unitlength}{\svgwidth}%
  \fi%
  \global\let\svgwidth\undefined%
  \global\let\svgscale\undefined%
  \makeatother%
  \begin{picture}(1,0.52912037)%
    \lineheight{1}%
    \setlength\tabcolsep{0pt}%
    \put(0,0){\includegraphics[width=\unitlength,page=1]{obtuse-1.pdf}}%
    \put(0.94628638,0.26745636){\color[rgb]{0,0,0}\makebox(0,0)[lt]{\lineheight{1.25}\smash{\begin{tabular}[t]{l}$x$\end{tabular}}}}%
    \put(0.78071802,0.35788212){\color[rgb]{0,0,0}\makebox(0,0)[lt]{\lineheight{1.25}\smash{\begin{tabular}[t]{l}$A$\end{tabular}}}}%
    \put(0.64291431,0.40653369){\color[rgb]{0,0,0}\makebox(0,0)[lt]{\lineheight{1.25}\smash{\begin{tabular}[t]{l}$a_2$\end{tabular}}}}%
    \put(0.64087651,0.27891882){\color[rgb]{0,0,0}\makebox(0,0)[lt]{\lineheight{1.25}\smash{\begin{tabular}[t]{l}$a_1$\end{tabular}}}}%
    \put(0.48014794,0.1209921){\color[rgb]{0,0,0}\makebox(0,0)[lt]{\lineheight{1.25}\smash{\begin{tabular}[t]{l}$B$\end{tabular}}}}%
    \put(0.44448718,0.42410945){\color[rgb]{0,0,0}\makebox(0,0)[lt]{\lineheight{1.25}\smash{\begin{tabular}[t]{l}$C'$\end{tabular}}}}%
    \put(0.18339871,0.50434643){\color[rgb]{0,0,0}\makebox(0,0)[lt]{\lineheight{1.25}\smash{\begin{tabular}[t]{l}$y$\end{tabular}}}}%
    \put(0.28401327,0.37189177){\color[rgb]{0,0,0}\makebox(0,0)[lt]{\lineheight{1.25}\smash{\begin{tabular}[t]{l}$y=f(x)$\end{tabular}}}}%
    \put(0,0){\includegraphics[width=\unitlength,page=2]{obtuse-1.pdf}}%
  \end{picture}%
\endgroup%

%% file: phi.pdf_tex
\begingroup%
  \makeatletter%
  \providecommand\color[2][]{%
    \errmessage{(Inkscape) Color is used for the text in Inkscape, but the package 'color.sty' is not loaded}%
    \renewcommand\color[2][]{}%
  }%
  \providecommand\transparent[1]{%
    \errmessage{(Inkscape) Transparency is used (non-zero) for the text in Inkscape, but the package 'transparent.sty' is not loaded}%
    \renewcommand\transparent[1]{}%
  }%
  \providecommand\rotatebox[2]{#2}%
  \newcommand*\fsize{\dimexpr\f@size pt\relax}%
  \newcommand*\lineheight[1]{\fontsize{\fsize}{#1\fsize}\selectfont}%
  \ifx\svgwidth\undefined%
    \setlength{\unitlength}{201.47109604bp}%
    \ifx\svgscale\undefined%
      \relax%
    \else%
      \setlength{\unitlength}{\unitlength * \real{\svgscale}}%
    \fi%
  \else%
    \setlength{\unitlength}{\svgwidth}%
  \fi%
  \global\let\svgwidth\undefined%
  \global\let\svgscale\undefined%
  \makeatother%
  \begin{picture}(1,0.59281201)%
    \lineheight{1}%
    \setlength\tabcolsep{0pt}%
    \put(0,0){\includegraphics[width=\unitlength,page=1]{phi.pdf}}%
    \put(0.42010592,0.07836679){\color[rgb]{0,0,0}\makebox(0,0)[lt]{\lineheight{1.25}\smash{\begin{tabular}[t]{l}$\delta$\end{tabular}}}}%
    \put(0.67217191,0.07935521){\color[rgb]{0,0,0}\makebox(0,0)[lt]{\lineheight{1.25}\smash{\begin{tabular}[t]{l}$2\delta$\end{tabular}}}}%
    \put(0.94104233,0.11988357){\color[rgb]{0,0,0}\makebox(0,0)[lt]{\lineheight{1.25}\smash{\begin{tabular}[t]{l}$x$\end{tabular}}}}%
    \put(0.13344259,0.56767141){\color[rgb]{0,0,0}\makebox(0,0)[lt]{\lineheight{1.25}\smash{\begin{tabular}[t]{l}$y$\end{tabular}}}}%
  \end{picture}%
\endgroup%